\documentclass[reqno,draft]{amsart}

\theoremstyle{plain}
\newtheorem{lemma}{Lemma}
\newtheorem{theorem}[lemma]{Theorem}

\newtheorem{definition}[lemma]{Definition}

\theoremstyle{remark}
\newtheorem{remark}{Remark}

\def\aa{\alpha}

\def\tt{\theta}

\def\Dd{\Delta}
\def\Gg{\Gamma}
\def\Om{\Omega}

\def\pp{\partial}

\begin{document}

\vskip 0.125in

\title[Primitive Equations]
{Global Well-posedness of the Three-dimensional
 Viscous Primitive Equations
of Large Scale Ocean and Atmosphere Dynamics}

\date{{\bf Submitted:} March 2, 2005. \quad {\bf Revised:}
 November 14, 2005}

\author[C. Cao]{Chongsheng Cao}
\address[C. Cao]
{Department of Mathematics  \\
Florida International University  \\
University Park  \\
Miami, FL 33199, USA.} \email{caoc@fiu.edu}

\author[E.S. Titi]{Edriss S. Titi}
\address[E.S. Titi]
{Department of Mathematics \\
and  Department of Mechanical and  Aerospace Engineering \\
University of California \\
Irvine, CA  92697-3875, USA. {\bf Also:}
 Department of Computer Science and Applied Mathematics \\
Weizmann Institute of Science  \\
Rehovot 76100, Israel.} \email{etiti@math.uci.edu}
\email{edriss.titi@weizmann.ac.il}

\begin{abstract}
In this paper we prove the global existence and uniqueness
(regularity)
 of strong solutions to the three-dimensional
 viscous primitive equations,
which model large scale  ocean and atmosphere  dynamics.
\end{abstract}

\maketitle


\vskip 0.125in 


\section{Introduction}   \label{S-1}

Large scale dynamics of oceans and atmosphere is governed by the
primitive equations which are derived from the
 Navier--Stokes equations, with rotation, coupled to
thermodynamics  and salinity diffusion-transport equations, which
account for the   buoyancy forces and
 stratification effects under the Boussinesq approximation.
Moreover, and  due to the shallowness of the oceans and the
atmosphere, i.e., the depth of  the fluid layer is very small in
comparison to the radius of the earth, the vertical large scale
motion in the oceans and the
 atmosphere is much smaller than the horizontal one,
which in turn leads to modeling the vertical motion by the
hydrostatic balance. As a result one obtains the
system~(\ref{EQ-1})--(\ref{EQ-4}), which is known as the primitive
equations for ocean and atmosphere dynamics (see, e.g.,
\cite{LTW92},\cite{LTW92A},\cite{PJ87}, \cite{Richardson},
\cite{S98}, \cite{TZ03}
 and references therein). We observe that  in the case of
ocean dynamics one has to add the diffusion-transport equation of
the salinity to the system~(\ref{EQ-1})--(\ref{EQ-4}). We omitted
it here in order
 to simplify our mathematical
presentation.  However, we  emphasize that our results are equally
valid when the salinity effects are taking into account.

 Let us remark that the horizontal motion can be further
approximated by the geostrophic balance when the Rossby number
(the ratio of the horizontal acceleration to the Coriolis force)
 is very small. By taking advantage of these assumptions and other
geophysical considerations several intermediate models
 have been developed and used in numerical studies
of weather prediction and  long-time climate dynamics (see, e.g.,
\cite{BKK}, \cite{CJ55A}, \cite{CJ55}, \cite{PJ87},
\cite{Richardson}, \cite{SR97}, \cite{SV97A}, \cite{SV97B},
\cite{SD96}, \cite{SH48} and references therein). Some of these
models have also been the subject of analytical  mathematical
study (see, e.g., \cite{BMN00}, \cite{BCT88}, \cite{CT03},
\cite{CTZ04}, \cite{CDGG}, \cite{CMT94A}, \cite{CMT94B},
\cite{EM98}, \cite{GMR01}, \cite{HTZ03},  \cite{J02},
 \cite{STW98}, \cite{STW00},
\cite{TZ03}, \cite{W88}  and references therein).

In this paper we will focus on  the $3D$ primitive equations
 in a cylindrical domain
\[
\Om= M\times (-h,0),
\]
where $M$ is a smooth bounded domain in $\mathbb{R}^2$:
\begin{eqnarray}
&&\hskip-.8in \frac{\pp v}{\pp t} + (v\cdot \nabla) v + w
\frac{\pp v}{\pp z} +
\nabla p + f \vec{k} \times v + L_1v = 0  \label{EQ-1}  \\
&&\hskip-.8in
\pp_z p  + T =0    \label{EQ-2}  \\
&&\hskip-.8in
\nabla \cdot v +\pp_z w =0   \label{EQ-3} \\
&&\hskip-.8in \frac{\pp T}{\pp t}  + v \cdot \nabla T + w
\frac{\pp T}{\pp z}  +  L_2 T =  Q \label{EQ-4}
\end{eqnarray}
where the horizontal velocity field $v=(v_1, v_2)$,  the
three-dimensional velocity field $(v_1, v_2, w)$, the temperature
$T$ and the pressure $p$ are the unknowns. $f=f_0(\beta+y)$ is the
Coriolis parameter, $Q$ is a given heat  source. The viscosity and
the heat diffusion operators $L_1$ and $L_2$ are given by
\begin{eqnarray}
&&\hskip-.8in
 L_1 = -\frac{1}{ Re_1} \Dd -  \frac{1}{  Re_2 } \;
\frac{\pp^2}{\pp z^2},   \label{L-1} \\
&&\hskip-.8in
 L_2 =  -\frac{1}{ Rt_1} \Dd -  \frac{1}{  Rt_2 } \;
\frac{\pp^2}{\pp z^2},  \label{L-2}
\end{eqnarray}
where $Re_1, Re_2$ are positive constants representing the
horizontal and vertical Reynolds numbers, respectively, and
$Rt_1, Rt_2$ are positive constants which stand for the horizontal
and vertical heat diffusivity,  respectively. We set $\nabla  =
(\pp_x, \pp_y)$ to be the horizontal gradient operator and $\Dd =
\pp_x^2 +\pp_y^2$ to be the horizontal Laplacian. We observe
that the above system is similar to the $3D$ Boussinesq 
system with the equation of vertical motion
is approximated by the hydrostatic balance. 

We partition the
boundary of $\Om$ into:
\begin{eqnarray}
&&\hskip-.8in
\Gg_u = \{ (x,y,z)  \in \overline{\Om} : z=0 \},  \\
&&\hskip-.8in
\Gg_b = \{ (x,y,z)  \in \overline{\Om} : z=-h \},  \\
&&\hskip-.8in \Gg_s = \{ (x,y,z)  \in \overline{\Om} : (x,y) \in
\partial M, \; -h \leq z \leq 0  \}.
\end{eqnarray}
We equip the system (\ref{EQ-1})--(\ref{EQ-4}) with the following
boundary conditions -- with wind--driven on the top surface and
non-slip and non-heat flux on the side walls and bottom (see,
e.g., {\bf \cite{LTW92}, \cite{LTW92A}, \cite{PJ87}, \cite{S98},
\cite{SR97}, \cite{SV97A},\cite{SV97B}, \cite{SD96}}):
\begin{eqnarray}
&&\hskip-.8in \mbox{on } \Gg_u: \frac{\pp v }{\pp z} = h\, \tau,
\; w=0, \;   \frac{\pp T }{\pp z} = - \aa (T-T^*);
\label{B-1}\\
&&\hskip-.8in \mbox{on } \Gg_b: \frac{\pp v }{\pp z} = 0, \; w=0,
\;
 \frac{\pp T }{\pp z} = 0;
\label{B-2}\\
&&\hskip-.8in \mbox{on } \Gg_s: v  \cdot \vec{n} =0,   \;
\frac{\pp v}{\pp \vec{n}}  \times \vec{n} =0, \; \frac{\pp T }{\pp
\vec{n}} =0,   \label{B-3}
\end{eqnarray}
where $\tau (x,y)$ is the wind stress on ocean surface, $\vec{n}$
is the normal vector to $\Gg_s$, and $T^* (x,y)$ is typical
temperature distribution of the top surface  of the ocean.
  For simplicity we assume here that
$\tau$ and $T^*$ are time independent. However, the results
presented here are equally valid when these quantities are time
dependent and satisfy certain
 bounds in space and time.

Due to the boundary conditions (\ref{B-1})--(\ref{B-3}), it is
natural to assume that $\tau$ and $T^*$ satisfy the compatibility
boundary  conditions:
\begin{eqnarray}
&&\hskip-.8in
 \tau \cdot \vec{n} =0,   \;  \frac{\pp \tau}{\pp \vec{n}}  \times
\vec{n} =0,   \qquad \mbox{on  } \pp M.
\label{COM-1}   \\
&&\hskip-.8in
 \frac{\pp T^* }{\pp \vec{n}} =0   \qquad
\mbox{on  } \pp M. \label{COM-2}
\end{eqnarray}

In addition, we supply the system with the initial condition:
\begin{eqnarray}
&&\hskip-.8in
v(x,y,z,0) = v_0 (x,y,z).  \label{INIT-1}\\
&&\hskip-.8in T(x,y,z,0) = T_0 (x,y,z).  \label{INIT-2}
\end{eqnarray}

In \cite{LTW92}, \cite{LTW92A} and \cite{TZ03} the authors set up
the mathematical framework to study the viscous primitive
equations for the atmosphere and ocean circulation. Moreover,
similar to the $3D$ Navier--Stokes equations, they have shown the
global existence of weak solutions, but the question of their
uniqueness is still open.  The short time existence and uniqueness
of strong solutions to the viscous primitive equations model was
established in \cite{GMR01} and \cite{TZ03}. In \cite{HTZ03} the
authors proved the global existence and uniqueness of strong
solutions to the viscous primitive equations in thin domains for a
large set of initial data whose size depends inversely on the
thickness of the domain. In this paper we show the global
existence, uniqueness and continuous dependence on initial
data, i.e.~global regularity and
well-posedness, of the strong
solutions to the $3D$ viscous primitive equations
model~(\ref{EQ-1})--(\ref{INIT-2}) in
 general cylindrical domain, $\Omega$,
and for any initial data. It is worth stressing that the
ideas developed in this paper can equally apply to
the primitive equations subject to other kids of boundary
conditions. As in the case of $3D$ Navier--Stokes
equations the question of uniqueness  of the weak solutions to
this model is still open.


\section{Preliminaries}    \label{S-2}

\vskip0.1in
\subsection{New Formulation} First, let us reformulate the system
(\ref{EQ-1})--(\ref{INIT-2}) (see also \cite{LTW92}, \cite{LTW92A}
and \cite{TZ03}). We integrate the equation (\ref{EQ-3}) in the
$z$ direction to obtain
\[
w(x,y,z,t) = w(x,y,-h,t) - \int_{-h}^z \nabla \cdot v(x,y, \xi,t)
d\xi.
\]
By virtue of (\ref{B-1}) and (\ref{B-2}) we have
\begin{equation}
w(x,y,z,t) =  - \int_{-h}^z \nabla \cdot v(x,y, \xi,t) d\xi,
\label{DIV-1}
\end{equation}
and
\begin{equation}
\int_{-h}^0 \nabla \cdot v(x,y, \xi,t) d\xi = \nabla \cdot
\int_{-h}^0 v(x,y, \xi,t) d\xi =0. \label{DIV}
\end{equation}
We denote by
\begin{equation}
\overline{\phi} (x, y)= \frac{1}{h} \int_{-h}^0 \phi(x,y, \xi)
d\xi, \qquad \forall \; (x,y) \in M. \label{VBAR}
\end{equation}
In particular,
\begin{equation}
\overline{v} (x, y)= \frac{1}{h} \int_{-h}^0  v(x,y, \xi) d\xi,
\qquad \mbox{in } \; M. \label{V--B}
\end{equation}
We will denote the fluctuation by
\begin{equation}
\widetilde{v} = v - \overline{v}.    \label{V--T}
\end{equation}
Notice that
\begin{eqnarray}
&&\hskip-.68in \overline{\widetilde{v}} = 0. \label{ZERO}
\end{eqnarray}
Based on the above and (\ref{B-3}) we obtain
\begin{eqnarray}
&&\hskip-.8in \nabla \cdot \overline{v} = 0,     \qquad \mbox{in }
\; M, \label{ADIV}
\end{eqnarray}
and
\begin{eqnarray}
&&\hskip-.8in \overline{v} \cdot \vec{n} = 0, \quad \frac{\pp
\overline{v}}{\pp \vec{n}} \times \vec{n} =0,  \qquad  \mbox{ on }
\pp M.   \label{AV-BD}
\end{eqnarray}
By integrating equation (\ref{EQ-2}) we obtain
\[
p(x,y,z,t) = - \int_{-h}^z T(x,y,\xi,t) d\xi + p_s(x,y,t).
\]
Substitute  (\ref{DIV-1}) and the above relation into equation
(\ref{EQ-1}) we reach
\begin{eqnarray}
&&\hskip-.68in \frac{\pp v}{\pp t} + (v\cdot \nabla) v - \left(
\int_{-h}^z \nabla \cdot v(x,y, \xi,t) d\xi
\right)  \frac{\pp v}{\pp  z}  \nonumber  \\
&&\hskip-.65in +  \nabla p_s(x,y,t) - \nabla  \int_{-h}^z
T(x,y,\xi,t) d\xi + f \vec{k} \times v + L_1 v = 0.   \label{VV}
\end{eqnarray}
\begin{remark}
Notice that due to the compatibility boundary  conditions
(\ref{COM-1}) and (\ref{COM-2}) one can convert the boundary
condition (\ref{B-1})--(\ref{B-3}) to be homogeneous by replacing
$(v, T)$ by $(v+\frac{(z+h)^2-h^3/3}{2} \tau, T+T^*)$ while
(\ref{ADIV}) is still true. For simplicity and without loss
generality we will assume that $\tau=0, T^*=0.$ However, we
emphasize that our results are still valid for general $\tau$ and
$T^*$ provided  they are smooth enough. In a forthcoming  paper we
will study the long-time dynamics and global attractors to the
primitive equations with general $\tau$ and $T^*$.
\end{remark}
Therefore, under the assumption that $\tau=0, T^*=0$, we have the
following new formulation for system (\ref{EQ-1})--(\ref{INIT-2}):
\begin{eqnarray}
&&\hskip-.68in \frac{\pp v}{\pp t} + L_1 v+ (v\cdot \nabla) v -
\left( \int_{-h}^z \nabla \cdot v(x,y, \xi,t) d\xi
\right)  \frac{\pp v}{\pp  z}  \nonumber  \\
&&\hskip-.65in +  \nabla p_s(x,y,t) - \nabla  \int_{-h}^z
T(x,y,\xi,t) d\xi
+ f \vec{k} \times v  = 0,   \label{EQV}   \\
&&\hskip-.68in
 \frac{\pp T}{\pp t}  + L_2 T   + v  \cdot \nabla T
- \left( \int_{-h}^z \nabla \cdot v (x,y, \xi,t)
 d\xi
\right) \frac{\pp T}{\pp z}   = Q,
 \label{EQ5}  \\
&&\hskip-.68in \left. \frac{\pp v }{\pp z} \right|_{z=0} = 0, \;
\left. \frac{\pp v }{\pp z} \right|_{z=-h} = 0, \; \left. v \cdot
\vec{n} \right|_{\Gg_s} = 0, \; \left. \frac{\pp v}{\pp \vec{n}}
\times \vec{n}\right|_{\Gg_s} =0,
\label{EQ6} \\
&&\hskip-.68in \left.  (\pp_z T + \aa T) \right|_{z=0}= 0; \;
\left. \pp_z T  \right|_{z=-h}= 0; \; \left. \pp_n T
\right|_{\Gg_s}= 0,
 \label{EQ7} \\
&&\hskip-.68in v (x,y,z,0) = v_0 (x,y,z),
\label{EQ8}   \\
&&\hskip-.68in T(x,y,z,0) = T_0 (x,y,z). \label{EQ9}
\end{eqnarray}

\subsection{Properties of $\overline{v}$ and $\widetilde{v}$}
By taking the average of equations (\ref{EQV}) in the $z$
direction, over the interval $(-h, 0)$, and using the boundary
conditions (\ref{EQ6}), we obtain
\begin{eqnarray}
&&\hskip-.68in \frac{\pp \overline{v}}{\pp t} + \overline{(v\cdot
\nabla) v - \left( \int_{-h}^z \nabla \cdot v(x,y, \xi,t) d\xi
\right) \frac{\pp v}{\pp   z}} + \nabla p_s(x,y,t) -\nabla  \left[
\frac{1}{h} \int_{-h}^0
\int_{-h}^z T(x,y,\xi,t) d\xi dz \right]   \nonumber  \\
&&\hskip-.65in + f \vec{k} \times \overline{v} - \frac{1}{Re_1}
\Dd \overline{v} = 0.   \label{VA}
\end{eqnarray}
As a result of (\ref{ZERO}), (\ref{ADIV}) and  integration by
parts we have
\begin{eqnarray}
&&\hskip-.268in
 \overline{(v\cdot \nabla) v - \left( \int_{-h}^z
\nabla \cdot v(x,y, \xi,t) d\xi \right) \frac{\pp v}{\pp   z}}=
(\overline{v} \cdot \nabla ) \overline{v} + \overline{
\left[(\widetilde{v} \cdot \nabla) \widetilde{v}
 + (\nabla \cdot \widetilde{v}) \; \widetilde{v}\right]}.
\label{I--1}
\end{eqnarray}
By subtracting  (\ref{VA}) from (\ref{EQV}) and using (\ref{I--1})
we get
\begin{eqnarray}
&&\hskip-.268in \frac{\pp \widetilde{v}}{\pp t} + L_1
\widetilde{v} + (\widetilde{v} \cdot \nabla) \widetilde{v} -
\left( \int_{-h}^z \nabla \cdot \widetilde{v}(x,y, \xi,t) d\xi
\right) \frac{\pp \widetilde{v}}{\pp  z} +(\widetilde{v} \cdot
\nabla ) \overline{v}+ (\overline{v} \cdot \nabla) \widetilde{v} -
\overline{ \left[(\widetilde{v} \cdot \nabla) \widetilde{v}  +
(\nabla
\cdot \widetilde{v}) \; \widetilde{v}\right]}  \nonumber  \\
&&\hskip-.265in
 - \nabla  \left( \int_{-h}^z T(x,y,\xi,t) d\xi
-\frac{1}{h} \int_{-h}^0 \int_{-h}^z T(x,y,\xi,t) d\xi dz  \right)
+ f \vec{k} \times \widetilde{v} =0.   \label{VTT}
\end{eqnarray}
Therefore,  $\overline{v}$ satisfies the following equations and
boundary conditions:
\begin{eqnarray}
&&\hskip-.68in \frac{\pp \overline{v}}{\pp t} - \frac{1}{Re_1} \Dd
\overline{v} + (\overline{v} \cdot \nabla ) \overline{v} +
\overline{ \left[ (\widetilde{v} \cdot \nabla) \widetilde{v}
 + (\nabla \cdot \widetilde{v}) \; \widetilde{v}\right]}
+ f \vec{k} \times \overline{v}   \nonumber  \\
&&\hskip-.6in + \nabla \left[ p_s(x,y,t) - \frac{1}{h} \int_{-h}^0
\int_{-h}^z  \, T (x,y,\xi,t) \;  d\xi \; dz \right]
= 0,  \label{EQ1}  \\
&&\hskip-.68in
\nabla \cdot \overline{v} = 0, \qquad \qquad
\mbox{ in } M,   \label{EQ2}     \\
&&\hskip-.68in \overline{v} \cdot \vec{n} = 0, \quad  \frac{\pp
\overline{v}}{\pp \vec{n}} \times \vec{n} =0,  \qquad \qquad
\mbox{ on } \pp M,   \label{EQ3}
\end{eqnarray}
and $\widetilde{v}$ satisfies the following equations and boundary
conditions:
\begin{eqnarray}
&&\hskip-.68in \frac{\pp \widetilde{v}}{\pp t} + L_1 \widetilde{v}
+ (\widetilde{v} \cdot \nabla) \widetilde{v} - \left( \int_{-h}^z
\nabla \cdot \widetilde{v}(x,y, \xi,t) d\xi \right) \frac{\pp
\widetilde{v}}{\pp  z} +(\widetilde{v} \cdot \nabla )
\overline{v}+ (\overline{v} \cdot \nabla) \widetilde{v}  \nonumber
\\
&&\hskip-.58in   - \overline{ \left[(\widetilde{v} \cdot \nabla)
\widetilde{v}  + (\nabla \cdot \widetilde{v}) \;
\widetilde{v}\right]} + f \vec{k} \times \widetilde{v} - \nabla
\left( \int_{-h}^z T(x,y,\xi,t) d\xi -\frac{1}{h} \int_{-h}^0
\int_{-h}^z
T(x,y,\xi,t) d\xi dz  \right)  =0,   \label{EQ4}    \\
&&\hskip-.68in \left. \frac{\pp \widetilde{v} }{\pp z}
\right|_{z=0} = 0, \; \left. \frac{\pp \widetilde{v} }{\pp z}
\right|_{z=-h} = 0, \; \left. \widetilde{v} \cdot \vec{n}
\right|_{\Gg_s} = 0, \; \left. \frac{\pp \widetilde{v}}{\pp
\vec{n}}  \times \vec{n}\right|_{\Gg_s} =0, \label{EQ66}.
\end{eqnarray}
\begin{remark} We recall that  by virtue of
the maximum principle one is able to show the global
well-posedness of the $3D$ viscous Burgers equations (see, for
instance, \cite{Lady} and references therein). Such an argument,
however, is not valid for the $3D$  Navier--Stokes equations
because of the pressure term. Remarkably, the pressure term is
absent from equation~(\ref{EQ4}). This fact  allows us
 to obtain a bound for the   $L^6$ norm of
$\widetilde{v}$, which is a key estimate in our proof of  the
global regularity for the
 system~(\ref{EQ-1})--(\ref{INIT-2}).

\end{remark}

\vskip0.1in
\subsection{Functional spaces and Inequalities}

Let us denote by $L^2(\Om), L^2(M)$ and $H^m(\Om), H^m(M)$ the
usual $L^2-$Lebesgue and Sobolev spaces, respectively
(\cite{AR75}). We denote by
\begin{equation}
\| \phi\|_p = \left\{ \begin{array}{ll} \left(  \int_{\Om}  |\phi
(x,y,z)|^p \; dxdydz \right)^{\frac{1}{p}},   \qquad & \mbox{ for
every $\phi \in L^p(\Om)$}
\\
\left(  \int_{M}  |\phi (x,y)|^p \; dxdy \right)^{\frac{1}{p}},  &
\mbox{ for every $\phi \in L^p(M)$}.
\end{array} \right.
 \label{L2}
\end{equation}
Let
\begin{eqnarray*}
\widetilde{\mathcal{V}_1} &=&  \left\{ v \in C^{\infty}(\Om):
\left. \frac{\pp v }{\pp z} \right|_{z=0} = 0, \; \left. \frac{\pp
v }{\pp z} \right|_{z=-h} = 0, \; \left. v \cdot \vec{n}
\right|_{\Gg_s} = 0, \; \left. \frac{\pp v}{\pp \vec{n}}  \times
\vec{n}\right|_{\Gg_s} =0,
\; \nabla \cdot \overline{v} =0 \right\},   \\
\widetilde{\mathcal{V}_2} &=&  \left\{ { T \in C^{\infty}(\Om):
\left. \frac{\pp T}{\pp z } \right|_{z=-h}= 0; \left. {\left(
\frac{\pp T}{\pp z} + \aa T \right) } \right|_{z= 0}= 0; \; \left.
\frac{\pp T}{\pp \vec{n}} \right|_{\Gg_s}= 0 } \right\}.
\end{eqnarray*}
We denote by $V_1$ and $V_2$ be the closure spaces of
$\widetilde{\mathcal{V}_1}$ in $H^1(\Om)$, and
$\widetilde{\mathcal{V}_2}$ in $H^1(\Om)$ under $H^1-$topology,
respectively.

\begin{definition} \label{D-1}
\thinspace Let  $v_0 \in V_1$ and $ T_0\in V_2$, and let
$\mathcal{T}$ be a fixed positive time. \thinspace $(v, T)$  is
called a strong solution of {\em (\ref{EQV})--(\ref{EQ9})} on the
time interval $[0,\mathcal{T}]$ if it satisfies {\em (\ref{EQV})}
and {\em (\ref{EQ5})} in weak sense, and also
\begin{eqnarray*}
&& v \in C([0,\mathcal{T}], V_1)  \cap L^2 ([0,\mathcal{T}],
H^2(\Om)),
\\
&&  T  \in C([0,\mathcal{T}], V_2)  \cap L^2 ([0,\mathcal{T}],
H^2(\Om)),    \\
&& \frac{dv}{dt} \in L^1([0,\mathcal{T}], L^2(\Om)),   \\
&& \frac{dT}{dt}  \; \in L^1([0,\mathcal{T}], L^2(\Om)).
\end{eqnarray*}

\end{definition}

For convenience, we recall the following Sobolev and
Ladyzhenskaya's inequalities in $\mathbb{R}^2$ (see, e.g.,
\cite{AR75}, \cite{CF88}, \cite{GA94}, \cite{LADY})
\begin{eqnarray}
&&\hskip-.68in \| \phi \|_{L^4(M)} \leq C_0 \| \phi \|_{L^2}^{1/2}
\| \phi
\|_{H^1(M)}^{1/2},  \label{SI-1}\\
&&\hskip-.68in \| \phi\|_{L^8(M)} \leq C_0 \| \phi
\|_{L^6(M)}^{3/4} \| \phi \|_{H^1(M)}^{1/4},  \label{SI-2}
\end{eqnarray}
for every $\phi \in H^1(M),$ and the following Sobolev and
Ladyzhenskaya's inequalities in $\mathbb{R}^3$ (see, e.g.,
\cite{AR75}, \cite{CF88}, \cite{GA94}, \cite{LADY})
\begin{eqnarray}
&&\hskip-.68in \| \psi \|_{L^3(\Om)} \leq C_0 \| u
\|_{L^2(\Om)}^{1/2} \| u
\|_{H^1(\Om)}^{1/2},  \label{SI1}\\
&&\hskip-.68in \| u \|_{L^6(\Om)} \leq C_0 \| u \|_{H^1(\Om)},
\label{SI2}
\end{eqnarray}
for every $u\in H^1(\Om).$  Here $C_0$ is a positive constant
which might depend on the shape of $M$ and $\Om$  but not on their
size. Moreover,  by (\ref{SI-1}) we get
\begin{eqnarray}
&&\hskip-.68in \| \phi \|_{L^{12}(M)}^{12} = \| |\phi|^3
\|_{L^{4}(M)}^{4} \leq C_0 \| |\phi|^3 \|_{L^2(M)}^{2} \| |\phi|^3
\|_{H^1(M)}^{2}
\nonumber  \\
&&\hskip-.68in \leq C_0  \| \phi \|_{L^6(M)}^{6} \; \left( \int_M
|\phi|^4 \left| \nabla \phi \right|^2 \; dxdy  \right) + \| \phi
\|_{L^6(M)}^{12},   \label{TWE}
\end{eqnarray}
for every $\phi \in H^1(M).$ Also, we recall the integral version
of Minkowsky inequality for the $L^p$ spaces, $p\geq 1$. Let
$\Om_1 \subset \mathbb{R}^{m_1}$ and
 $\Om_2 \subset \mathbb{R}^{m_2}$ be two measurable sets, where
$m_1$ and $m_2$ are two positive integers. Suppose that
$f(\xi,\eta)$ is measurable over $\Om_1 \times \Om_2$. Then,
\begin{equation}
\hskip0.35in \left[ { \int_{\Om_1} \left( \int_{\Om_2}
|f(\xi,\eta)| d\eta \right)^p d\xi } \right]^{1/p} \leq
\int_{\Om_2} \left( \int_{\Om_1} |f(\xi,\eta)|^p d\xi
\right)^{1/p} d\eta. \label{MKY}
\end{equation}

\section{{\em A Priori} estimates} \label{S-3}

In the previous subsections we have reformulated the system
(\ref{EQ-1})--(\ref{INIT-2}) and obtained the  system
(\ref{EQV})--(\ref{EQ9}). The two systems are equivalent when
$(v,T)$ is a strong solution. The existence of such a strong
solution for a short interval of time, whose length depends on the
initial data and the other physical parameters of the system
(\ref{EQ-1})--(\ref{INIT-2}), was established in \cite{ GMR01} and
\cite{TZ03}. Let $(v_0,T_0)$ be a given initial data. In this
section we will consider the strong solution that corresponds  to
this initial data in its maximal interval of existence
$[0,\mathcal{T_*})$. Specifically, we will establish {\em a
priori} upper estimates for various norms of this solution in the
interval $[0,\mathcal{T_*})$. In particular, we will show that if
$\mathcal{T_*} < \infty$ then the $H^1$ norm of the strong
solution is bounded over the interval $[0,\mathcal{T_*})$. This
key observation  plays a major role in the proof of global
regularity of strong solutions to the system
(\ref{EQ-1})--(\ref{INIT-2}).

\subsection{$L^2$ estimates}

We take  the inner product of equation  (\ref{EQ5}) with $T$, in
$L^2(\Om)$,  and obtain
\begin{eqnarray*}
&&\hskip-.68in \frac{1}{2} \frac{d \|T\|_2^2}{dt} + \frac{1}{Rt_1}
\|\nabla T\|_2^2
+ \frac{1}{Rt_2}\|T_z\|_2^2  +\aa \|T(z=0)\|_2^2 \\
&&\hskip-.65in = \int_{\Om} Q T \; dxdydz  -\int_{\Om} \left( v
\cdot \nabla T - \left( \int_{-h}^z \nabla \cdot v(x,y, \xi,t)
d\xi \right) \frac{\pp T}{\pp z}\right) T \; dxdydz.
\end{eqnarray*}
After integrating by parts we get
\begin{eqnarray}
&&\hskip-.065in   -\int_{\Om} \left( v \cdot \nabla T - \left(
\int_{-h}^z \nabla \cdot v(x,y, \xi,t) d\xi \right) \frac{\pp
T}{\pp z}\right) T \; dxdydz =0. \label{DT-9}
\end{eqnarray}
As a result of the above we conclude
\begin{eqnarray*}
&&\hskip-.68in \frac{1}{2} \frac{d \|T\|_2^2}{dt} + \frac{1}{Rt_1}
\|\nabla T\|_2^2 + \frac{1}{Rt_2}\|T_z\|_2^2  +\aa \|T(z=0)\|_2^2 \\
&&\hskip-.65in =\int_{\Om} Q T \; dxdydz \leq \|Q\|_2 \; \|T\|_2.
\end{eqnarray*}
Notice that
\begin{equation}
\| T \|_2^2 \leq 2 h^2 \| T_z \|_2^2 + 2 h \|T (z=0)\|_2^2.
\label{P2}
\end{equation}
Using (\ref{P2}) and  the Cauchy--Schwarz inequality we obtain
\begin{eqnarray}
&&\hskip-.68in  \frac{d \|T\|_2^2}{dt} + \frac{2}{Rt_1}
\|\nabla T\|_2^2 + \frac{1}{Rt_2}\|T_z\|_2^2  +\aa \|T(z=0)\|_2^2 \\
&&\hskip-.65in  \leq  2(h^2\; Rt_2 + \frac{h}{\aa}) \|Q\|^2_2.
\label{T_E}
\end{eqnarray}
By the inequality (\ref{P2}) and thanks to Gronwall inequality the
above gives
\begin{eqnarray}
&&\hskip-.68in  \|T\|_2^2  \leq
  e^{-\; \frac{t}{2(h^2\; Rt_2 + h/\aa)}} \|T_0\|_2^2 +
(2 h^2\; Rt_2 + 2h/\aa)^2 \|Q\|^2_2,   \label{T-2}
\end{eqnarray}
Moreover, we have
\begin{eqnarray}
&&\hskip-.68in   \int_0^t \left[  \frac{1}{Rt_1} \|\nabla T
(s)\|_2^2 + \frac{1}{Rt_2}\|T_z(s)\|^2_2  +\aa \|T(z=0)(s)\|^2_2
\right]\; ds
\nonumber \\
&&\hskip-.65in \leq  2(h^2\; Rt_2 + \frac{h}{\aa}) \|Q\|^2_2 \; t
+
 e^{-\; \frac{t}{2(h^2\; Rt_2 + h/\aa)}} \|T_0\|_2^2 +
(2 h^2\; Rt_2 + 2h/\aa)^2 \|Q\|^2_2,   \label{T-2I}
\end{eqnarray}
By taking the inner product of equation (\ref{EQV}) with $v$, in
$L^2(\Om)$,  we reach
\begin{eqnarray*}
&&\hskip-.28in \frac{1}{2} \frac{d \|v\|_2^2}{dt} + \frac{1}{Re_1}
\|\nabla v\|_2^2
+ \frac{1}{Re_2}\|v_z\|_2^2   \\
&&\hskip-.265in =  -\int_{\Om}  \left[  (v \cdot \nabla) v -
\left( \int_{-h}^z \nabla \cdot v(x,y, \xi,t) d\xi \right)
\frac{\pp v}{\pp  z} \right] \cdot  v \;dxdydz
 \\
 &&\hskip-.1065in +\int_{\Om}  \left(
f \vec{k} \times v +\nabla p_s - \nabla  \left( \int_{-h}^z
T(x,y,\xi,t) d\xi  \right) \right) \cdot v \;dxdydz.
\end{eqnarray*}
By integration by parts we get
\begin{eqnarray}
&&\hskip-.065in \int_{\Om}  \left[ (v \cdot \nabla) v - \left(
\int_{-h}^z \nabla \cdot v(x,y,\xi,t) d\xi \right) \frac{\pp
v}{\pp  z} \right] \cdot v \;dxdydz =0.   \label{DT-1}
\end{eqnarray}
By (\ref{EQ2}) we have
\begin{eqnarray}
&&\hskip-.065in \int_{\Om} \nabla p_s \cdot v \;dxdydz = h \int_M
\nabla p_s \cdot \overline{v} \;dxdy = -h \int_{\Om}  p_s (\nabla
\cdot \overline{v}) \;dxdy =0. \label{DT-4}
\end{eqnarray}
Since
\begin{eqnarray}
&&\hskip-.065in  (f \vec{k} \times v )\cdot v=0, \label{DT-6}
\end{eqnarray}
then from (\ref{DT-1})--(\ref{DT-6}) we have
\begin{eqnarray*}
&&\hskip-.68in \frac{1}{2} \frac{d \|v \|_2^2 }{dt} +
\frac{1}{Re_1} \|\nabla v\|_2^2
+ \frac{1}{Re_2}\|v_z\|_2^2    \\
&&\hskip-.65in =   - \int_{\Om}  \int_{-h}^z T(x,y,\xi,t)\; d\xi
(\nabla  \cdot v ) \; dxdydz   \\
&&\hskip-.68in \leq  h \|T\|_2\; \| \nabla v \|_2.
\end{eqnarray*}
 By Cauchy--Schwarz and (\ref{T-2}) we obtain
\begin{eqnarray*}
&&\hskip-.68in   \frac{1}{2} \frac{d \|v \|_2^2 }{dt} +
\frac{1}{Re_1} \|\nabla v\|_2^2
+ \frac{1}{Re_2}\|v_z\|_2^2    \\
&&\hskip-.68in \leq  h^2 Re_1 \; \|T\|_2^2  \leq h^2 Re_1 \left(
\|T_0\|_2^2 + (2 h^2\; Rt_2 + 2h/\aa)^2 \|Q\|^2_2\right).
\end{eqnarray*}
Recall that (cf., e.g., \cite{GA94} Vol. I  p. 55)
\[
\| v \|^2_2 \leq C_M \| \nabla v \|^2_2.
\]
By the above and thanks to Gronwall inequality we get
\begin{eqnarray}
&&\hskip-.68in \| v\|_2^2 \leq e^{- \; \frac{t}{C_M \; Re_1 } }
\left( h \| \overline{v_0}\|_2^2 + \| \widetilde{v_0}\|_2^2
\right)  \nonumber   \\
&&\hskip-.568in  + C_M h^2 Re_1^2 \left[ \;  \|T_0\|_2^2 + (2h^2\;
Rt_2 + 2h/\aa)^2 \|Q\|^2_2\; \right]. \label{VEE}
\end{eqnarray}
Moreover,
\begin{eqnarray}
&&\hskip-.68in \int_0^t \left[ \frac{1}{Re_1} \| \nabla v(s)
\|_2^2   + \frac{1}{Re_2} \|  v_z (s)\|_2^2  \right]\; ds
\nonumber \\
&&\hskip-.65in \leq h^2 Re_1 \left( \|T_0\|_2^2 + (2 h^2\; Rt_2 +
2h/\aa)^2 \|Q\|^2_2\right) \; t + e^{- \; \frac{t}{C_M \; Re_1 } }
\left( h \| \overline{v_0}\|_2^2 + \| \widetilde{v_0}\|_2^2
\right)  \nonumber   \\
&&\hskip-.568in  + C_M h^2 Re_1^2 \left[\; \|T_0\|_2^2 + (2h^2\;
Rt_2 + 2h/\aa)^2 \|Q\|^2_2\; \right]. \label{VEEI}
\end{eqnarray}
Therefore, by (\ref{T-2}), (\ref{T-2I}), (\ref{VEE}) and
(\ref{VEEI}) we have
\begin{eqnarray}
&&\hskip-.168in  \| v (t)\|_2^2 + \int_0^t \left[ \frac{1}{Re_1}
\| \nabla v (s)\|_2^2  +
\frac{1}{Re_2} \|  v_z (s)\|_2^2  \right]\; ds   \nonumber   \\
&&\hskip-.068in + \|T (t)\|_2^2 + \int_0^t \left[  \frac{1}{Rt_1}
\|\nabla T (s)\|_2^2 + \frac{1}{Rt_2}\|T_z(s)\|^2_2  +\aa
\|T(z=0)(s)\|^2_2  \right]\; ds   \leq K_1 (t),   \label{K-1}
\end{eqnarray}
where
\begin{eqnarray}
&&\hskip-.668in K_1 (t) = 2(h^2\; Rt_2 + h/\aa) \|Q\|^2_2 ; t  +
\left( h \| \overline{v_0}\|_2^2 + \| \widetilde{v_0}\|_2^2
\right)  \nonumber   \\
&&\hskip-.568in  + \left( 1+ C_M h^2 Re_1^2 + h^2 Re_1 \, t
\right) \left[\; \|T_0\|_2^2 + (2h^2\; Rt_2 + 2h/\aa)^2
\|Q\|^2_2\; \right]. \label{K1}
\end{eqnarray}

\subsection{$L^6$ estimates}

Taking the  inner product of the  equation (\ref{EQ4}) with
$|\widetilde{v}|^4 \widetilde{v}$ in $L^2(\Om)$,  we get
\begin{eqnarray*}
&&\hskip-.168in \frac{1}{6} \frac{d \| \widetilde{v} \|_{6}^{6}
}{d t} + \frac{1}{Re_1} \int_{\Om} \left(|\nabla \widetilde{v}|^2
|\widetilde{v}|^{4} + \left|\nabla |\widetilde{v}|^2 \right|^2
|\widetilde{v}|^{2} \right) \; dxdydz + \frac{1}{Re_2} \int_{\Om}
\left(|\widetilde{v}_z|^2 |\widetilde{v}_z|^{4} + \left|\pp_z
|\widetilde{v}|^2 \right|^2 |\widetilde{v}|^{2} \right)
\; dxdydz    \\
&&\hskip-.165in = - \int_{\Om} \left\{  (\widetilde{v} \cdot
\nabla) \widetilde{v} - \left( \int_{-h}^z \nabla \cdot
\widetilde{v}(x,y, \xi,t) d\xi \right)  \frac{\pp
\widetilde{v}}{\pp  z} +(\widetilde{v} \cdot \nabla )
\overline{v}+ (\overline{v} \cdot \nabla) \widetilde{v} -
\overline{ \left[(\widetilde{v} \cdot \nabla) \widetilde{v}
+ (\nabla \cdot \widetilde{v}) \; \widetilde{v}\right]} \right.  \\
&&\hskip-.06in
 \left.
+ f \vec{k} \times \widetilde{v} - \nabla  \left( \int_{-h}^z
T(x,y,\xi,t) d\xi -\frac{1}{h} \int_{-h}^0 \int_{-h}^z
T(x,y,\xi,t) d\xi dz  \right)    \right\} \cdot
|\widetilde{v}|^{4} \widetilde{v} \; dxdydz.
\end{eqnarray*}
By integration by parts we get
\begin{eqnarray}
&&\hskip-.065in   - \int_{\Om} \left[ (\widetilde{v} \cdot \nabla)
\widetilde{v} - \left( \int_{-h}^z \nabla \cdot \widetilde{v}(x,y,
\xi,t) d\xi \right)  \frac{\pp \widetilde{v}}{\pp  z} \right]
\cdot |\widetilde{v}|^{4} \widetilde{v} \; dxdydz =0. \label{D6-1}
\end{eqnarray}
Since
\begin{eqnarray}
&&\hskip-.065in  \left(f \vec{k} \times \widetilde{v} \right)
\cdot |\widetilde{v}|^{4} \widetilde{v}=0, \label{D6-2}
\end{eqnarray}
then by (\ref{EQ2}) and the boundary condition (\ref{EQ6}) we also
have
\begin{eqnarray}
&&\hskip-.065in \int_{\Om} (\overline{v} \cdot \nabla)
\widetilde{v}\cdot |\widetilde{v}|^{4} \widetilde{v} \;dxdydz =0.
\label{D6-3}
\end{eqnarray}
Thus, by (\ref{D6-1})--(\ref{D6-3}) we have
\begin{eqnarray*}
&&\hskip-.168in \frac{1}{6} \frac{d \| \widetilde{v} \|_{6}^{6}
}{d t} + \frac{1}{Re_1} \int_{\Om} \left(|\nabla \widetilde{v}|^2
|\widetilde{v}|^{4} + \left|\nabla |\widetilde{v}|^2 \right|^2
|\widetilde{v}|^{2} \right) \; dxdydz + \frac{1}{Re_2} \int_{\Om}
\left(|\widetilde{v}_z|^2 |\widetilde{v}_z|^{4} + \left|\pp_z
|\widetilde{v}|^2 \right|^2 |\widetilde{v}|^{2} \right)
\; dxdydz    \\
&&\hskip-.165in
 = - \int_{\Om} \left\{
(\widetilde{v} \cdot \nabla ) \overline{v}  - \overline{
(\widetilde{v} \cdot \nabla) \widetilde{v}
+ (\nabla \cdot \widetilde{v}) \; \widetilde{v}}   \right.  \\
&&\hskip-.06in
 \left.
- \nabla  \left( \int_{-h}^z T(x,y,\xi,t) d\xi -\frac{1}{h}
\int_{-h}^0 \int_{-h}^z T(x,y,\xi,t) d\xi dz \right) \right\}
\cdot |\widetilde{v}|^{4} \widetilde{v} \; dxdydz.
\end{eqnarray*}
Notice that by integration by parts and boundary condition
(\ref{EQ6}) we have
\begin{eqnarray*}
&&\hskip-.165in
 - \int_{\Om} \left[
(\widetilde{v} \cdot \nabla ) \overline{v}  - \overline{\left[
(\widetilde{v} \cdot \nabla) \widetilde{v}
+ (\nabla \cdot \widetilde{v}) \; \widetilde{v} \right]}   \right.  \\
&&\hskip-.06in \left. - \nabla  \left( \int_{-h}^z T(x,y,\xi,t)
d\xi -\frac{1}{h} \int_{-h}^0 \int_{-h}^z T(x,y,\xi,t) d\xi dz
\right)    \right] \cdot
|\widetilde{v}|^{4} \widetilde{v} \; dxdydz  \\
&&\hskip-.165in
 = \int_{\Om} \left[
(\nabla \cdot \widetilde{v})  \; \overline{v} \cdot
|\widetilde{v}|^{4} \widetilde{v} + (\widetilde{v} \cdot \nabla )
(|\widetilde{v}|^{4} \widetilde{v}) \cdot \overline{v}
 - \overline{ \widetilde{v}^k   \widetilde{v}^j} \;  \pp_{x_k}
(|\widetilde{v}|^{4} \widetilde{v}^j) \right.  \\
&&\hskip-.06in
 \left.
-  \left( \int_{-h}^z T(x,y,\xi,t) d\xi -\frac{1}{h} \int_{-h}^0
\int_{-h}^z T(x,y,\xi,t) d\xi dz  \right) \nabla \cdot
(|\widetilde{v}|^{4} \widetilde{v})   \right]  \; dxdydz.
\end{eqnarray*}
Therefore, by Cauchy--Schwarz inequality and H\"{o}lder inequality
we obtain
\begin{eqnarray*}
&&\hskip-.168in \frac{1}{6} \frac{d \| \widetilde{v} \|_{6}^{6}
}{d t} + \frac{1}{Re_1} \int_{\Om} \left(|\nabla \widetilde{v}|^2
|\widetilde{v}|^{4} + \left|\nabla |\widetilde{v}|^2 \right|^2
|\widetilde{v}|^{2} \right) \; dxdydz + \frac{1}{Re_2} \int_{\Om}
\left(|\widetilde{v}_z|^2 |\widetilde{v}_z|^{4} + \left|\pp_z
|\widetilde{v}|^2 \right|^2 |\widetilde{v}|^{2} \right)
\; dxdydz    \\
&&\hskip-.165in \leq C \int_{M} \left[ |\overline{v}| \int_{-h}^0
|\nabla \widetilde{v}| \; | \widetilde{v}|^5 \; dz \right]
\; dxdy  \\
&&\hskip-.06in +C \int_{M} \left[ \left( \int_{-h}^0
|\widetilde{v}|^2 \; dz \right) \left( \int_{-h}^0 |\nabla
\widetilde{v}| \;
| \widetilde{v}|^4 \; dz \right) \right] \; dxdy \\
&&\hskip-.065in +C \int_{M} \left[ \overline{|T|} \int_{-h}^0
|\nabla \widetilde{v}| \; | \widetilde{v}|^4 \; dz \right]
\; dxdy   \\
&&\hskip-.165in \leq C \int_{M} \left[ |\overline{v}| \left(
\int_{-h}^0 |\nabla \widetilde{v}|^2 \; | \widetilde{v}|^4 \; dz
\right)^{1/2} \left( \int_{-h}^0  | \widetilde{v}|^6 \; dz
\right)^{1/2} \right]
\; dxdy  \\
&&\hskip-.06in +C \int_{M} \left[ \left( \int_{-h}^0
|\widetilde{v}|^2 \; dz \right) \left( \int_{-h}^0 |\nabla
\widetilde{v}|^2 \; | \widetilde{v}|^4 \; dz \right)^{1/2} \left(
\int_{-h}^0  | \widetilde{v}|^4 \; dz \right)^{1/2} \right] \;
dxdy \\
&&\hskip-.06in +C \int_{M} \left[ \overline{|T|} \left(
\int_{-h}^0 |\nabla \widetilde{v}|^2 \; | \widetilde{v}|^4 \; dz
\right)^{1/2} \left( \int_{-h}^0  | \widetilde{v}|^4 \; dz
\right)^{1/2} \right]
\; dxdy   \\
&&\hskip-.165in \leq C \|\overline{v}\|_{L^4(M)}
 \left( \int_{\Om}  |\nabla \widetilde{v}|^2 \; | \widetilde{v}|^4 \;
dxdydz \right)^{1/2} \left( \int_M \left( \int_{-h}^0  |
\widetilde{v}|^6 \; dz \right)^{2}
\; dxdy \right)^{1/4}   \\
&&\hskip-.06in +C \left( \int_M \left( \int_{-h}^0  |
\widetilde{v}|^2 \; dz \right)^{4} \; dxdy \right)^{1/4} \left(
\int_{\Om}  |\nabla \widetilde{v}|^2 \; | \widetilde{v}|^4 \;
dxdydz \right)^{1/2} \left( \int_M \left( \int_{-h}^0  |
\widetilde{v}|^4 \; dz \right)^{2}
\; dxdy \right)^{1/4} \\
&&\hskip-.06in +C \|\,\overline{|T|}\,\|_{L^4(M)} \left(
\int_{\Om} |\nabla \widetilde{v}|^2 \; | \widetilde{v}|^4 \;
dxdydz \right)^{1/2} \left( \int_M \left( \int_{-h}^0  |
\widetilde{v}|^4 \; dz \right)^{2} \; dxdy \right)^{1/4}.
\end{eqnarray*}
By using Minkowsky  inequality (\ref{MKY}),  we get
\begin{eqnarray*}
&&\hskip-.168in \left( \int_M \left( \int_{-h}^0  |
\widetilde{v}|^6 \; dz \right)^{2} \; dxdy \right)^{1/2} \leq  C
\int_{-h}^0 \left( \int_M  | \widetilde{v}|^{12} \; dxdy
\right)^{1/2}  \; dz.
\end{eqnarray*}
By (\ref{TWE}),
\begin{eqnarray*}
&&\hskip-.168in \int_M  | \widetilde{v}|^{12}  \; dxdy \leq C_0
\left( \int_M | \widetilde{v}|^{6}  \; dxdy \right) \left( \int_M
| \widetilde{v}|^{4} |\nabla \widetilde{v} |^2  \; dxdy \right) +
\left( \int_M | \widetilde{v}|^{6}  \; dxdy \right)^2.
\end{eqnarray*}
Thus, by Cauchy--Schwarz inequality we obtain
\begin{eqnarray}
&&\hskip-.168in \left( \int_M \left( \int_{-h}^0  |
\widetilde{v}|^6 \; dz \right)^{2} \; dxdy \right)^{1/2} \leq  C
\| \widetilde{v}\|_{L^6(\Om)}^{3}  \left( \int_{\Om}
\widetilde{v}|^{4} |\nabla \widetilde{v} |^2 \; dxdydz
\right)^{1/2} + \| \widetilde{v}\|_{L^6(\Om)}^{6}. \label{M1}
\end{eqnarray}
Similarly, by (\ref{MKY}) and (\ref{SI-2}), we also get
\begin{eqnarray}
&&\hskip-.168in \left( \int_M \left( \int_{-h}^0  |
\widetilde{v}|^4 \; dz \right)^{2} \; dxdy \right)^{1/2} \leq  C
\int_{-h}^0 \left( \int_M  | \widetilde{v}|^{8}
\; dxdy \right)^{1/2}  \; dz   \nonumber  \\
&&\hskip-.168in \leq  C  \int_{-h}^0  \|
\widetilde{v}\|_{L^6(M)}^{3} \left( \|\nabla
\widetilde{v}\|_{L^2(M)}  + \| \widetilde{v}\|_{L^2(M)} \right) \;
dz \leq C  \| \widetilde{v}\|_6^{3} \left( \|\nabla
\widetilde{v}\|_2 + \| \widetilde{v}\|_2 \right), \label{M2}
\end{eqnarray}
and
\begin{eqnarray}
&&\hskip-.168in \left( \int_M \left( \int_{-h}^0  |
\widetilde{v}|^2 \; dz \right)^{4} \; dxdy \right)^{1/4} \leq  C
\int_{-h}^0 \left( \int_M  | \widetilde{v}|^{8}
\; dxdy \right)^{1/4}  \; dz    \nonumber  \\
&&\hskip-.168in \leq  C  \int_{-h}^0  \|
\widetilde{v}\|_{L^6(M)}^{3/2}\left(
 \|\nabla \widetilde{v}\|_{L^2(M)}^{1/2}   + \|
\widetilde{v}\|_{L^2(M)}^{1/2} \right) \; dz \leq C  \|
\widetilde{v}\|_6^{3/2} \left( \|\nabla \widetilde{v}\|_2^{1/2}  +
\| \widetilde{v}\|_2^{1/2} \right). \label{M3}
\end{eqnarray}
Therefore, by (\ref{M1})--(\ref{M3}) and (\ref{SI-1}), we reach
\begin{eqnarray*}
&&\hskip-.168in \frac{1}{6} \frac{d \| \widetilde{v} \|_{6}^{6}
}{d t} + \frac{1}{Re_1} \int_{\Om} \left(|\nabla \widetilde{v}|^2
|\widetilde{v}|^{4} + \left|\nabla |\widetilde{v}|^2 \right|^2
|\widetilde{v}|^{2} \right) \; dxdydz + \frac{1}{Re_2} \int_{\Om}
\left(|\widetilde{v}_z|^2 |\widetilde{v}_z|^{4} + \left|\pp_z
|\widetilde{v}|^2 \right|^2 |\widetilde{v}|^{2} \right)
\; dxdydz    \\
&&\hskip-.165in \leq C \|\overline{v}\|_2^{1/2} \; \|\nabla
\overline{v}\|_2^{1/2} \| \widetilde{v}\|_6^{3/2}
 \left( \int_{\Om}  |\nabla \widetilde{v}|^2 \; | \widetilde{v}|^4 \;
dxdydz \right)^{3/4} +C \|\overline{v}\|_2^{1/2} \; \|\nabla
\overline{v}\|_2^{1/2}  \| \widetilde{v}\|_6^{6}
   \\
&&\hskip-.065in +C \| \widetilde{v}\|_6^{3} \left( \|\nabla
\widetilde{v}\|_2 + \| \widetilde{v}\|_2 \right) \left( \int_{\Om}
|\nabla \widetilde{v}|^2 \; | \widetilde{v}|^4 \;
dxdydz \right)^{1/2}   \\
&&\hskip-.065in +C \|\overline{T}\|_2^{1/2} \; \|\nabla
\overline{T}\|_2^{1/2} \| \widetilde{v}\|_6^{3/2} \left( \|\nabla
\widetilde{v}\|_2^{1/2}  + \| \widetilde{v}\|_2^{1/2} \right)
\left( \int_{\Om}  |\nabla \widetilde{v}|^2 \; | \widetilde{v}|^4
\; dxdydz \right)^{1/2}.
\end{eqnarray*}
Thanks to the Young's and the Cauchy--Schwarz inequalities we have
\begin{eqnarray*}
&&\hskip-.168in
 \frac{d \| \widetilde{v} \|_{6}^{6} }{d t} +
\frac{1}{Re_1} \int_{\Om} \left(|\nabla \widetilde{v}|^2
|\widetilde{v}|^{4} + \left|\nabla |\widetilde{v}|^2 \right|^2
|\widetilde{v}|^{2} \right) \; dxdydz + \frac{1}{Re_2} \int_{\Om}
\left(|\widetilde{v}_z|^2 |\widetilde{v}_z|^{4} + \left|\pp_z
|\widetilde{v}|^2 \right|^2 |\widetilde{v}|^{2} \right)
\; dxdydz    \\
&&\hskip-.165in \leq C \|\overline{v}\|_2^{2} \; \|\nabla
\overline{v}\|_2^{2}   \| \widetilde{v}\|_6^{6} +C \|
\widetilde{v}\|_6^{6} \|\nabla \widetilde{v}\|_2^2 +C
\|\overline{T}\|_2^2 \; \|\nabla \overline{T}\|_2^2 +C \|
\widetilde{v}\|_2^2\| \widetilde{v}\|_6^{6}.
\end{eqnarray*}
By (\ref{K-1}) and Gronwall inequality, we get
\begin{eqnarray}
&&\hskip-.68in \| \widetilde{v} (t)\|^6_6 + \int_0^t \left(
\frac{1}{Re_1} \int_{\Om} |\nabla \widetilde{v}|^2
|\widetilde{v}|^{4}
 \; dxdydz +
\frac{1}{Re_2} \int_{\Om} |\widetilde{v}_z|^2
|\widetilde{v}_z|^{4}  \; dxdydz \right) \leq K_6 (t), \label{K-6}
\end{eqnarray}
where
\begin{eqnarray}
&&\hskip-.68in K_6 (t) = e^{K_1^2 (t)} \left[ \|v_0\|_{H^1(\Om)}^6
+ K_1^2 (t) \right].   \label{K6}
\end{eqnarray}

Taking the  inner product of the  equation (\ref{EQ5}) with $|T|^4
T$ in $L^2(\Om)$, and by (\ref{EQ5}),  we get
\begin{eqnarray*}
&&\hskip-.168in \frac{1}{6} \frac{d \| T \|_{6}^{6} }{d t} +
\frac{5}{Rt_1} \int_{\Om} |\nabla T|^2 |T|^{4}  \; dxdydz +
\frac{5}{Rt_2} \int_{\Om} |T_z|^2 |T|^{4} \; dxdydz  +\aa
\|T(z=0)\|_6^6 \\
&&\hskip-.165in = \int_{\Om} Q |T|^4 T \; dxdydz  -\int_{\Om}
\left( v
 \cdot \nabla T
- \left( \int_{-h}^z \nabla \cdot v(x,y, \xi,t)
 d\xi
\right) \frac{\pp T}{\pp z}\right) |T|^4 T \; dxdydz.
\end{eqnarray*}
By integration by parts and (\ref{EQ2}) we get
\begin{eqnarray}
&&\hskip-.065in   -\int_{\Om} \left( v \cdot \nabla T - \left(
\int_{-h}^z \nabla \cdot v(x,y, \xi,t)
 d\xi
\right) \frac{\pp T}{\pp z}\right) |T|^4 T \; dxdydz =0.
\label{DT-10}
\end{eqnarray}
As a result of the above we conclude
\begin{eqnarray*}
&&\hskip-.168in \frac{1}{6} \frac{d \| T \|_{6}^{6} }{d t} +
\frac{5}{Rt_1} \int_{\Om} |\nabla T|^2 |T|^{4}  \; dxdydz +
\frac{5}{Rt_2} \int_{\Om} |T_z|^2 |T|^{4} \; dxdydz  +\aa
\|T(z=0)\|_6^6  \\
&&\hskip-.165in = \int_{\Om} Q|T|^4 T\; dxdydz \leq \|Q\|_6
\|T\|_6^5.
\end{eqnarray*}
By Gronwall, again, we get
\begin{eqnarray}
&&\hskip-.68in \| T  (t)\|_6  \leq \|Q\|_{H^1(\Om)} \; t + \|
T_0\|_{H^1(\Om)}.     \label{K-T}
\end{eqnarray}

\subsection{$H^1$ estimates}

\subsubsection{$\|\nabla \overline{v}\|_2$ estimates}

First, let us observe that since $v$ is a strong solution
on the interval $[0,\mathcal{T}_*)$  then
$\Dd\overline{v}\in  L^2 ([0,\mathcal{T}_*), L^2(M))$. 
Consequently, and  by virtue of~(\ref{EQ2}), 
$\Dd\overline{v}\cdot \vec{n}\in  L^2 ([0,\mathcal{T}_*), 
H^{-1/2}(\pp M))$ (see, e.g., \cite{CF88}, \cite{TT84}).
Moreover, and thanks to~(\ref{EQ2}) and~(\ref{EQ3}), we
have  $\Dd\overline{v}\cdot \vec{n}=0$ on $\pp M$ (see, e.g.,
\cite{Ziane98}). This observation implies also that the
Stokes operator in the domain $M$, subject to the boundary
conditions~(\ref{EQ3}), is equal to the $-\Dd$ operator.

As a result of the above and~(\ref{EQ2})
we apply a generalized version of the Stokes theorem 
(see, e.g.,  \cite{CF88}, \cite{TT84}) to conclude:
\[
\int_M \nabla p_s(x,y,t)\cdot \Dd\overline{v} (x,y,t)dx dy =0.
\]

By taking the inner product of equation~(\ref{EQ1}) with $- \Dd
\overline{v}$ in $L^2(M)$, and applying~(\ref{EQ2}) and the
above, we reach
\begin{eqnarray*}
&&\hskip-.68in \frac{1}{2} \frac{d \| \nabla \overline{v} \|_2^2
}{d t} + \frac{1}{Re_1} \|\Dd \overline{v}\|_2^2  =  \int_{M}
\left\{ (\overline{v} \cdot \nabla ) \overline{v} + \overline{
\left[ (\widetilde{v} \cdot \nabla) \widetilde{v}  + (\nabla \cdot
\widetilde{v}) \; \widetilde{v}\right]} \right\} \cdot \Dd
\overline{v} \; dxdy  + \int_{M} f \vec{k} \times \overline{v}
\cdot \Dd \overline{v} \; dxdy.
\end{eqnarray*}
Following similar steps as in the proof 
of $2D$ Navier--Stokes
equations (cf. e.g., \cite{CF88}, \cite{TT84}) one obtains
\begin{eqnarray*}
&&\hskip-.68in \left| \int_{M}  (\overline{v} \cdot \nabla )
\overline{v} \cdot \Dd \overline{v} \; dxdy \right|
 \leq C \|\overline{v}\|_2^{1/2} \|\nabla \overline{v}\|_2
\; \|\Dd  \overline{v}\|_2^{3/2}.
\end{eqnarray*}
Applying the Cauchy--Schwarz  and H\"{o}lder inequalities, we get
\begin{eqnarray*}
&&\hskip-.68in \left|\int_{M} \; \overline{ (\widetilde{v} \cdot
\nabla) \widetilde{v}  + (\nabla \cdot \widetilde{v}) \;
\widetilde{v}}  \cdot \Dd \overline{v} \; dxdy \right| \leq C
\int_M \int_{-h}^0 |\widetilde{v}| \; |\nabla \widetilde{v}| \; dz
\; |\Dd \overline{v}| \; dxdy
\\
&&\hskip-.68in \leq C \int_M \left[ \left( \int_{-h}^0
|\widetilde{v}|^2 \; |\nabla \widetilde{v}| \; dz \right)^{1/2}
\left( \int_{-h}^0  |\nabla \widetilde{v}| \; dz \right)^{1/2} \;
|\Dd \overline{v}| \right] \; dxdy
\\
&&\hskip-.68in \leq C \left[ \int_M  \left( \int_{-h}^0
|\widetilde{v}|^2 \; |\nabla \widetilde{v}| \; dz \right)^{2} \;
dxdy \right]^{1/4} \; \left[ \int_M \left( \int_{-h}^0 |\nabla
\widetilde{v}| \; dz \right)^{2} \; dxdy \right]^{1/4} \; \left[
\int_M
 |\Dd \overline{v}|^2  \; dxdy  \right]^{1/2} \\
&&\hskip-.68in \leq C \|\nabla \widetilde{v}\|_2^{1/2} \left(
\int_{\Om} |\widetilde{v}|^4 |\nabla \widetilde{v}|^2 \; dxdydz
\right)^{1/4}\|\Dd \overline{v}\|_2.
\end{eqnarray*}
Thus, by Young's and Cauchy--Schwarz inequalities, we have
\begin{eqnarray*}
&&\hskip-.68in \frac{d \| \nabla \overline{v} \|_2^2 }{d t} +
\frac{1}{Re_1} \|\Dd \overline{v}\|_2^2 \leq C
\|\overline{v}\|_2^2 \|\nabla \overline{v}\|_2^4 + C \|\nabla
\widetilde{v}\|_2^{2} + C \int_{\Om} |\widetilde{v}|^4 |\nabla
\widetilde{v}|^2 \; dxdydz + C \|\overline{v}\|_2^2.
\end{eqnarray*}
By (\ref{K-1}), (\ref{K-6}) and thanks to Gronwall inequality and
we obtain
\begin{eqnarray}
&&\hskip-.68in \| \nabla \overline{v} \|_2^2 + \frac{1}{Re_1}
\int_0^t |\Dd \overline{v}|_2^2 \; ds \leq K_2(t),   \label{K-2}
\end{eqnarray}
where
\begin{eqnarray}
&&\hskip-.68in K_2 (t)= e^{K_1^2 (t)} \left[ \|v_0\|_{H^1(\Om)}^2
+ K_1(t) +K_6(t) \right].   \label{K2}
\end{eqnarray}

\subsubsection{$\|v_z\|_2$ estimates}

Denote by $u=v_z.$ It is clear that $u$ satisfies
\begin{eqnarray}
&&\hskip-.68in
 \frac{\pp u }{\pp t} + L_1 u +
(v \cdot \nabla) u - \left( \int_{-h}^z \nabla \cdot v (x,y,
\xi,t)
d\xi \right) \frac{\pp u}{\pp z}    \nonumber    \\
&&\hskip-.65in + (u \cdot \nabla ) v - (\nabla \cdot v) u +f
\vec{k} \times u -\nabla T
  = 0. \label{UU}
\end{eqnarray}
Taking the  inner product of the  equation (\ref{UU}) with $u$ in
$L^2$ and using the boundary condition (\ref{EQ6}),  we get
\begin{eqnarray*}
&&\hskip-.68in \frac{1}{2} \frac{d \| u \|_2^2 }{d t} +
\frac{1}{Re_1} \|\nabla u\|_2^2 +
\frac{1}{Re_2}  \|\pp_z u\|_2^2     \\
&&\hskip-.65in =- \int_{\Om} \left(  (v \cdot \nabla) u - \left(
\int_{-h}^z \nabla \cdot v (x,y, \xi,t)
d\xi \right) \frac{\pp u}{\pp z}  \right) \cdot u \; dxdydz \\
&&\hskip-.58in - \int_{\Om} \left( (u \cdot \nabla ) v - (\nabla
\cdot v) u +f \vec{k} \times u -\nabla T \right) \cdot u \;
dxdydz.
\end{eqnarray*}
By integration by parts we get
\begin{eqnarray}
&&\hskip-.065in   - \int_{\Om} \left( (v \cdot \nabla) u - \left(
\int_{-h}^z \nabla \cdot v(x,y, \xi,t) d\xi \right)  \frac{\pp
u}{\pp  z} \right) \cdot u \; dxdydz =0. \label{DZ-1}
\end{eqnarray}
Since
\begin{eqnarray}
&&\hskip-.065in  (f \vec{k} \times u) \cdot u=0, \label{DZ-2}
\end{eqnarray}
then by (\ref{DZ-1}) and (\ref{DZ-2}) we have
\begin{eqnarray*}
&&\hskip-.68in \frac{1}{2} \frac{d \| u \|_2^2 }{d t} +
\frac{1}{Re_1} \|\nabla u\|_2^2 +
\frac{1}{Re_2}  \|\pp_z u\|_2^2     \\
&&\hskip-.65in =- \int_{\Om} \left(  (u \cdot \nabla ) v - (\nabla
\cdot v) u -\nabla
T  \right) \cdot u \; dxdydz \\
&&\hskip-.65in \leq C \int_{\Om} \left( |v| \right) \,|u| \,
|\nabla u| \; dxdydz
+ \|T\|_2 \; \|\nabla u\|_2    \\
&&\hskip-.65in \leq C  \|v\|_6 \; \| u \|_3  \|\nabla
u\|_2 + \|T\|_2\; \|\nabla u\|_2  \\
&&\hskip-.65in \leq C \|v\|_6\; \| u \|_2^{1/2}  \|\nabla
u\|_2^{3/2} + \|T\|_2\; \|\nabla u\|_2.
\end{eqnarray*}
By Young's inequality and Cauchy--Schwarz inequality, we have
\begin{eqnarray*}
&&\hskip-.68in \frac{d \| u \|_2^2 }{d t} + \frac{1}{Re_1}
\|\nabla u\|_2^2 +
\frac{1}{Re_2}  \|\pp_z u\|_2^2     \\
&&\hskip-.65in  \leq C \|v\|_6^4 \; \| u \|_2^{2}  + C \|T\|_2^2   \\
&&\hskip-.65in  \leq C \left( \|\nabla \overline{v}\|_2^4
+\|\widetilde{v}\|_6^4\right) \; \| u \|_2^{2}  + C \|T\|_2^2.
\end{eqnarray*}
By  (\ref{K-1}), (\ref{K-6}), (\ref{K-2}), and Gronwall
inequality,  we get
\begin{eqnarray}
&&\hskip-.68in \| v_z \|_2^2 + \frac{1}{Re_1} \int_0^t \|\nabla
v_z (s)\|_2^2 +\frac{1}{Re_2} \int_0^t \|v_{zz} (s)\|_2^2 \; ds
\leq K_z (t),   \label{K-Z}
\end{eqnarray}
where
\begin{eqnarray}
&&\hskip-.68in K_z (t)= e^{(K_2^2(t)+K_6^{2/3}(t)) t} \left[
\|v_0\|_{H^1(\Om)}^2 + K_1 (t)\right].   \label{KZ}
\end{eqnarray}

\subsubsection{$\|\nabla v\|_2$ estimates}

By taking the inner product of equation (\ref{EQV}) with $- \Dd v$
in $L^2(\Om)$, we reach
\begin{eqnarray*}
&&\hskip-.68in \frac{1}{2} \frac{d \| \nabla v \|_2^2}{d t} +
\frac{1}{Re_1} \|\Dd v\|_2^2
+\frac{1}{Re_2} \|\nabla v_z\|_2^2 \\
&&\hskip-.68in = - \int_{\Om} \left[(v \cdot \nabla) v - \left(
\int_{-h}^z \nabla \cdot v(x,y, \xi,t) d\xi \right) \frac{\pp
v}{\pp  z}   \right.    \\
&&\hskip-.528in \left. + f \vec{k} \times v +\nabla p_s- \nabla
\left( \int_{-h}^z
T(x,y,\xi,t) d\xi  \right)  \right] \cdot \Dd v \; dxdydz  \\
&&\hskip-.65in \leq C \int_{\Om} \left[ |v| \,|\nabla v|
 + \int_{-h}^0 |\nabla v| \; dz |\widetilde{v}_z|
+
 \int_{-h}^0 |\nabla T|\; dz \; \right] \, |\Dd v| \; dxdydz  \\
&&\hskip-.65in \leq C \|v\|_{L^6(\Om)}\; \| \nabla v\|_{L^3(\Om)}
\|\Dd v\|_2
 + C \int_M \left( \int_{-h}^0 |\nabla  v|
\; dz \int_{-h}^0 |v_z| |\Dd v| \; dz \right) \; dxdy + C \|\nabla
T\|_2 \; \|\Dd v\|_2.
\end{eqnarray*}
Notice that by applying the Proposition 2.2 in {\bf \cite{CT03}}
with $u=v, f=\Dd v$ and $g=v_z$, we get
\begin{eqnarray*}
&&\hskip-.68in \int_M \left( \int_{-h}^0 |\nabla  v| \; dz
\int_{-h}^0 |v_z| |\Dd v| \; dz \right) \; dxdy \leq C \|\nabla
v\|_2^{1/2} \|v_z\|_2^{1/2} \|\nabla v_z\|_2^{1/2} \|\Dd
v\|_2^{3/2}.
\end{eqnarray*}
As a result and by (\ref{SI1}) and (\ref{SI2}), we obtain
\begin{eqnarray*}
&&\hskip-.68in \frac{1}{2} \frac{d \| \nabla v \|_2^2 }{d t} +
\frac{1}{Re_1} \|\Dd v\|_2^2
+\frac{1}{Re_2} \|\nabla v_z\|_2^2 \\
&&\hskip-.65in \leq C \left( \|v \|_{L^6(\Om)} +\|\nabla
v\|_2^{1/2} \|v_z\|_2^{1/2} \right) \; \| \nabla v \|_2^{1/2}
\|\Dd v\|_2^{3/2}
 + h \|\nabla T\|_2 \; \|\Dd v\|_2.
\end{eqnarray*}
Thus, by Young's inequality and Cauchy--Schwarz inequality, we
have
\begin{eqnarray*}
&&\hskip-.68in \frac{d \| \nabla v \|_2^2 }{d t} + \frac{1}{Re_1}
\|\Dd v\|_2^2
+\frac{1}{Re_2} \|\nabla v_z\|_2^2 \\
&&\hskip-.68in \leq C \left( \|v\|_{L^6(\Om)}^4 +\|\nabla
v\|_2^{2} \|v_z\|_2^{2}  \right) \|\nabla v \|_2^2
 + C \|\nabla T\|_2^2.
\end{eqnarray*}
By (\ref{K-1}), (\ref{K-6}), (\ref{K-2}), (\ref{K-Z}) and thanks
to Gronwall inequality,  we obtain
\begin{eqnarray}
&&\hskip-.68in \| \nabla  v\|_2^2 + \int_0^t \left( \frac{1}{Re_1}
\|\Dd v(s)\|_2^2 +\frac{1}{Re_2} \|\nabla v_z(s)\|_2^2 \right) \;
ds \leq K_V(t), \label{K-V}
\end{eqnarray}
where
\begin{eqnarray}
&&\hskip-.68in K_V (t)= e^{K_6^{2/3}(t) \; t + K_1(t)\; K_z(t) }
\left[ \|v_0\|_{H^1(\Om)}^2 + K_1 (t) \right].   \label{KV}
\end{eqnarray}

\subsubsection{$\|T\|_{H^1}$ estimates}

Taking the  inner product of the  equation (\ref{EQ5}) with $-\Dd
T -T_{zz}$ in $L^2(\Om)$,  we get
\begin{eqnarray*}
&&\hskip-.168in \frac{1}{2} \frac{d \left( \|\nabla T\|_2^2 +\|
T_z \|_2^2
+ \aa \|\nabla T(z=0)\|_2^2\right) }{d t} \\
&&\hskip-.016in  + \frac{1}{Rt_1} \|\nabla T_z\|_2^2 + \left(
\frac{1}{Rt_1} +\frac{1}{Rt_2} \right) \left( \|\nabla T_z \|_2^2
+\aa \|\nabla T(z=0)\|_2^2 \right)+
\frac{1}{Rt_2}  \|T_{zz}\|_2^2     \\
&&\hskip-.165in = \int_{\Om} \left[  v \cdot \nabla T
 - \left(\int_{-h}^z \nabla \cdot v \; d\xi \right) T_z  -Q \right]
\; \left[ \Dd T+ T_{zz} \right]\; dxdydz  \\
&&\hskip-.165in \leq C \int_{\Om} \left(  |v| \, |\nabla T|   +|Q|
\right) \; \left|\Dd T+ T_{zz} \right|\; dxdydz + \int_M \left[
\int_{-h}^0 |\nabla v| \; dz \; \int_{-h}^0 |T_z| \; \left|\Dd T+
T_{zz} \right|\; dz \right]
\; dxdy \\
&&\hskip-.165in \leq C \|v\|_6 \; \| \nabla T \|_3  \left( \|\Dd
T\|_2^2+ \|\nabla T_z\|_2^2
+ \|T_{zz}\|_2^2 \right)^{1/2} \\
&&\hskip-.0165in +  C \|\nabla v\|_2^{1/2} \|\Dd v\|_2^{1/2}
\|T_z\|_2^{1/2} \left( \|\Dd T\|_2^2+ \|\nabla T_z\|_2^2+
\|T_{zz}\|_2^2 \right)^{3/2} +\|Q\|_2 \left( \|\Dd T\|_2^2+
\|\nabla T_z\|_2^2+ \|T_{zz}\|_2^2
\right)^{1/2}  \\
&&\hskip-.165in \leq C \left[ \|v\|_6  \; \| \nabla T \|_2^{1/2}
+ \|\nabla v\|_2^{1/2} \|\Dd v\|_2^{1/2} \|T_z\|_2^{1/2} \right]
\left( \|\Dd T\|_2^2+ \|\nabla T_z\|_2^2+ \|T_{zz}\|_2^2
\right)^{3/2}
\\
&&\hskip-.165in +\|Q\|_2 \left( \|\Dd T\|_2^2+ \|\nabla T_z\|_2^2+
\|T_{zz}\|_2^2 \right)^{1/2}.
\end{eqnarray*}
By Young's inequality and Cauchy--Schwarz inequality  we have
\begin{eqnarray*}
&&\hskip-.68in \frac{d \left(\|\nabla T\|_2^2 +\| T_z \|_2^2 + \aa
\|\nabla T(z=0)\|_2^2\right) }{d t} \\
&&\hskip-.6in  + \frac{1}{Rt_1} \|\nabla T_z\|_2^2 + \left(
\frac{1}{Rt_1} +\frac{1}{Rt_2} \right) \left( \|\nabla T_z \|_2^2
+\aa \|\nabla T(z=0)\|_2^2 \right) +
\frac{1}{Rt_2}  \|T_{zz}\|_2^2     \\
&&\hskip-.65in \leq C \left( \|v\|_6^4  +\|\nabla v\|_2^{2} \|\Dd
v\|_2^{2}  \right) \; \left( \|\nabla T\|_2^2 +\| T_z
\|_2^2\right) +C \|Q\|_2^2.
\end{eqnarray*}
By (\ref{K-6}), (\ref{K-V}), and  Gronwall inequality, we get
\begin{eqnarray}
&&\hskip-.68in \|\nabla T\|_2^2 +\| T_z \|_2^2 + \aa \|\nabla
T(z=0)\|_2^2   \nonumber
\\
&&\hskip-.6in +  \int_0^t\left[ \frac{1}{Rt_1} \|\nabla T_z\|_2^2
+ \left( \frac{1}{Rt_1} +\frac{1}{Rt_2} \right) \left( \|\nabla
T_z \|_2^2 +\aa \|\nabla T(z=0)\|_2^2 \right) \frac{1}{Rt_2}
\|T_{zz}\|_2^2 \right] \; ds   \leq K_t,   \label{V-T}
\end{eqnarray}
where
\begin{eqnarray}
&&\hskip-.68in K_t = e^{K_6^2(t)\, t + K_V^2(t)} \left[
\|T_0\|_{H^1(\Om)}^2 + \|Q\|_2^2 \right].   \label{VT}
\end{eqnarray}

\section{Existence and Uniqueness of the Strong Solutions} \label{S-4}

In this section  we will use the  {\em a priori} estimates
(\ref{K-1})--(\ref{V-T}) to show the global existence and
uniqueness, i.e. global regularity, of strong solutions to the
system (\ref{EQV})--(\ref{EQ9}).

\begin{theorem} \label{T-MAIN}
Let $Q \in H^1(\Om)$,  $v_0 \in V_1$, $ T_0\in V_2$ and
$\mathcal{T}>0,$ be given. Then there exists  a unique strong
solution $(v, T)$ of the system {\em (\ref{EQV})--(\ref{EQ9})} on
the interval $[0,\mathcal{T}]$ which depends continuously on the
initial data.
\end{theorem}

\vskip0.05in

\begin{proof}
As we have indicated earlier the short time existence of the
strong solution was established in \cite{GMR01} and \cite{TZ03}.
Let $(v,T)$ be the strong solution corresponding to the initial
data $(v_0,T_0)$ with maximal interval of existence
$[0,\mathcal{T}_*)$. If we assume that $\mathcal{T_*} < \infty$
then it is clear that
\[
\limsup_{t \to \mathcal{T}_*^{-}} \left( \| v \|_{H^1(\Om)}+\| T
\|_{H^1(\Om)} \right) =\infty.
\]
Otherwise, the solution can be extended beyond the time
$\mathcal{T_*}$. However, the above contradicts the {\it a priori}
estimates (\ref{K-Z}), (\ref{K-V}) and (\ref{V-T}). Therefore
$\mathcal{T_*}= \infty$, and the solution $(v, T)$ exists globally
in time.

Next,  we show the continuous dependence on the initial data and
the  the uniqueness of the strong solutions. Let $(v_1, T_1)$ and
$(v_2,T_2)$ be two strong solutions of the system
(\ref{EQV})--(\ref{EQ9}) with corresponding pressures $(p_s)_1$
and  $(p_s)_2$, and initial data $((v_0)_1, (T_0)_1)$ and
$((v_0)_2, (T_0)_2)$, respectively. Denote by $u=v_1-v_2, q_s
=(p_s)_1 -(p_s)_2$ and $ \tt = T_1-T_2.$ It is clear that
\begin{eqnarray}
&&\hskip-.68in \frac{\pp u}{\pp t} + L_1 u + (v_1 \cdot \nabla) u
+ (u \cdot \nabla) v_2 - \left( \int_{-h}^z \nabla \cdot v_1(x,y,
\xi,t) d\xi \right) \frac{\pp u}{\pp  z} - \left( \int_{-h}^z
\nabla \cdot u(x,y, \xi,t) d\xi \right) \frac{\pp v_2}{\pp  z}
\nonumber  \\
&&\hskip-.6in
 + f \vec{k} \times u + \nabla q_s - \nabla  \left( \int_{-h}^z
\tt (x,y,\xi,t) d\xi   \right)  =0,   \label{UEQ1}    \\
&&\hskip-.68in
 \frac{\pp \tt}{\pp t}  + L_2 \tt  + v_1
 \cdot \nabla \tt + u
 \cdot \nabla T_2   \nonumber  \\
&&\hskip-.6in - \left( \int_{-h}^z \nabla \cdot v_1(x,y, \xi,t)
 d\xi
\right) \frac{\pp \tt }{\pp z} - \left( \int_{-h}^z \nabla \cdot
u(x,y, \xi,t)
 d\xi
\right) \frac{\pp T_2}{\pp z}   = 0,
 \label{UEQ2}   \\
&&\hskip-.68in u(x,y,z,t) = (v_0)_1-(v_0)_2,  \label{UEQ3}   \\
&&\hskip-.68in \tt (x,y,z,0) = (T_0)_1-(T_0)_2. \label{UEQ4}
\end{eqnarray}
By taking the inner product of equation (\ref{UEQ1}) with with $u$
in $L^2(\Om)$, and equation (\ref{UEQ2}) with $\tt$, in $L^2(\Om)$
we get
\begin{eqnarray*}
&&\hskip-.268in \frac{1}{2} \frac{d \|u\|_2^2}{dt}
+ \frac{1}{Re_1} \|\nabla u\|_2^2 + \frac{1}{Re_2}\|u_z\|_2^2  \\
&&\hskip-.265in =   - \int_{\Om} \left[ (v_1 \cdot \nabla) u + (u
\cdot \nabla) v_2 - \left( \int_{-h}^z \nabla \cdot v_1(x,y,
\xi,t) d\xi \right) \frac{\pp u}{\pp  z} - \left( \int_{-h}^z
\nabla \cdot u(x,y, \xi,t) d\xi \right) \frac{\pp v_2}{\pp  z}
\right]
\cdot u \; dxdydz    \\
&&\hskip-.16in
 - \int_{\Om} \left[  f \vec{k} \times u + \nabla q_s - \nabla  \left(
\int_{-h}^z \tt (x,y,\xi,t) d\xi   \right) \right] \cdot u \;
dxdydz,
\end{eqnarray*}
and
\begin{eqnarray*}
&&\hskip-.168in \frac{1}{2} \frac{d \|\tt\|_2^2}{dt} +
\frac{1}{Rt_1} \|\nabla \tt\|_2^2 + \frac{1}{Rt_2}\|\tt_z\|_2^2
+\aa
\|\tt (z=0)\|_2^2 \\
&&\hskip-.165in =   - \int_{\Om} \left[ v_1
 \cdot \nabla \tt + u
 \cdot \nabla T_2
- \left( \int_{-h}^z \nabla \cdot v_1(x,y, \xi,t)
 d\xi
\right) \frac{\pp \tt }{\pp z} - \left( \int_{-h}^z \nabla \cdot
u(x,y, \xi,t)
 d\xi
\right) \frac{\pp T_2}{\pp z}\right]\; \tt \; dxdydz.
\end{eqnarray*}
By integration by parts, and the boundary conditions (\ref{EQ6})
and (\ref{EQ7}), we get
\begin{eqnarray}
&&\hskip-.065in   - \int_{\Om} \left( (v_1 \cdot \nabla) u -
\left( \int_{-h}^z \nabla \cdot v_1(x,y, \xi,t) d\xi \right)
\frac{\pp u}{\pp  z} \right) \cdot u \; dxdydz =0, \label{DUU-1}   \\
&&\hskip-.065in   - \int_{\Om} \left( v_1 \cdot \nabla \tt
 - \left( \int_{-h}^z \nabla \cdot v_1(x,y, \xi,t) d\xi \right)
\frac{\pp \tt}{\pp  z} \right) \cdot \tt \; dxdydz =0.
\label{DUT-1}
\end{eqnarray}
Since
\begin{eqnarray}
&&\hskip-.065in  \left(f \vec{k} \times u\right) \cdot u=0,
\label{DUU-2}
\end{eqnarray}
Then by (\ref{DUU-1}), (\ref{DUT-1}) and (\ref{DUU-2}) we have
\begin{eqnarray*}
&&\hskip-.68in \frac{1}{2} \frac{d \|u\|_2^2}{dt}
+ \frac{1}{Re_1} \|\nabla u\|_2^2 + \frac{1}{Re_2}\|u_z\|_2^2  \\
&&\hskip-.65in =   - \int_{\Om} (u \cdot \nabla) v_2  \cdot u \;
dxdydz + \int_{\Om}  \int_{-h}^z \nabla \cdot u(x,y, \xi,t) d\xi
\frac{\pp v_2}{\pp  z} \cdot u \; dxdydz.
\end{eqnarray*}
and
\begin{eqnarray*}
&&\hskip-.68in \frac{1}{2} \frac{d \|\tt\|_2^2}{dt} +
\frac{1}{Rt_1} \|\nabla \tt\|_2^2 + \frac{1}{Rt_2}\|\tt_z\|_2^2
+\aa
\|\tt (z=0)\|_2^2 \\
&&\hskip-.65in =   - \int_{\Om} (u \cdot \nabla) T_2  \tt \;
dxdydz + \int_{\Om}  \int_{-h}^z \nabla \cdot u(x,y, \xi,t) d\xi
\frac{\pp T_2}{\pp  z} \tt \; dxdydz.
\end{eqnarray*}
Notice that
\begin{eqnarray}
&&\hskip-.68in \left| \int_{\Om} (u \cdot \nabla) v_2  \cdot u \;
dxdydz \right| \leq \|\nabla v_2\|_2 \|u\|_3 \|u\|_6 \leq C
\|\nabla v_2\|_2 \|u\|_2^{1/2} \|\nabla u\|_2^{3/2},
\label{U1}   \\
&&\hskip-.68in \left| \int_{\Om} (u \cdot \nabla) T_2  \tt \;
dxdydz \right| \leq \|\nabla v_2\|_2 \|\tt \|_3 \|u\|_6 \leq C
\|\nabla T_2\|_2 \|\tt\|_2^{1/2} \|\nabla \tt\|_2^{1/2} \|\nabla
u\|_2.   \label{U-1}
\end{eqnarray}
Moreover,
\begin{eqnarray*}
&&\hskip-.38in \left| \int_{\Om}  \int_{-h}^z \nabla \cdot u(x,y,
\xi,t) d\xi \frac{\pp v_2}{\pp  z} \cdot u \; dxdydz \right| \leq
\int_M  \left( \int_{-h}^0 |\nabla u|\,  dz  \int_{-h}^0 |\pp_z
v_2| \, |u| \,  dz \right) \; dxdy \\
&&\hskip-.35in \leq  \int_M  \left( \int_{-h}^0 |\nabla u|\,  dz
\left( \int_{-h}^0 |\pp_z v_2|^2  \,  dz\right)^{1/2}
 \left( \int_{-h}^0  |u|^2 \,  dz\right)^{1/2}  \right) \; dxdy   \\
&&\hskip-.35in \leq  \left( \int_M  \left( \int_{-h}^0 |\nabla
u|\, dz \right)^2 \; dxdy \right)^{\frac{1}{2}} \left( \int_M
\left( \int_{-h}^0 |\pp_z v_2|^2  \,  dz\right)^{2} \; dxdy
\right)^{\frac{1}{4}} \left( \int_M  \left( \int_{-h}^0  |u|^2 \,
dz\right)^{2}   \; dxdy \right)^{\frac{1}{4}}.
\end{eqnarray*}
By Cauchy--Schwarz inequality, we get
\begin{eqnarray}
&&\hskip-.168in \left( \int_M  \left( \int_{-h}^0 |\nabla u|\,  dz
\right)^2 \; dxdy \right)^{1/2} \leq C \|\nabla u\|_2. \label{U-2}
\end{eqnarray}
By using Minkowsky  inequality (\ref{MKY}) and (\ref{SI-1}), we
obtain
\begin{eqnarray}
&&\hskip-.168in \left( \int_M \left( \int_{-h}^0  | u|^2 \; dz
\right)^{2} \; dxdy \right)^{1/2} \leq  C  \int_{-h}^0 \left(
\int_M  | u|^{4}
\; dxdy \right)^{1/2}  \; dz    \nonumber   \\
&&\hskip-.168in \leq C  \int_{-h}^0 | u| |\nabla u|  \; dz   \leq
C \| u\|_2 \|\nabla u\|_2, \label{U2}
\end{eqnarray}
and
\begin{eqnarray}
&&\hskip-.168in \left( \int_M \left( \int_{-h}^0  | \pp_z v_2|^2
\; dz \right)^{2} \; dxdy \right)^{1/2} \leq  C  \int_{-h}^0
\left( \int_M  | \pp_z v_2|^{4}
\; dxdy \right)^{1/2}  \; dz    \nonumber   \\
&&\hskip-.168in \leq C  \int_{-h}^0 | \pp_z v_2| |\nabla \pp_z
v_2|  \; dz   \leq C \| \pp_z v_2\|_2 \|\nabla \pp_z v_2\|_2.
\label{U3}
\end{eqnarray}
Similarly, we have
\begin{eqnarray}
&&\hskip-.68in \left| \int_{\Om}  \int_{-h}^z \nabla \cdot u(x,y,
\xi,t) d\xi \frac{\pp T_2}{\pp  z} \tt \; dxdydz \right| \leq C
\|\nabla u\|_2 \| \pp_z T_2\|_2^{1/2} \|\nabla \pp_z T_2\|_2^{1/2}
\|\tt\|_2^{1/2} \|\nabla \tt\|_2^{1/2}. \label{U-3}
\end{eqnarray}
Therefore, by estimates (\ref{U1})--(\ref{U-3}), we reach
\begin{eqnarray*}
&&\hskip-.68in \frac{1}{2} \frac{d
\left(\|u\|_2^2+\|\tt\|_2^2\right) }{dt} + \frac{1}{Re_1} \|\nabla
u\|_2^2 + \frac{1}{Re_2}\|u_z\|_2^2 + \frac{1}{Rt_1} \|\nabla
\tt\|_2^2 + \frac{1}{Rt_2}\|\tt_z\|_2^2 +\aa
\|\tt (z=0)\|_2^2\\
&&\hskip-.65in \leq C \left( \|\nabla v_2\|_2+ \| \pp_z
v_2\|_2^{1/2} \|\nabla \pp_z v_2\|_2^{1/2}  \right)
\|u\|^{1/2} \|\nabla u\|_2^{3/2} \\
&&\hskip-.6in +   C\|\nabla T_2\|_2 \|\tt\|_2^{1/2} \|\nabla
\tt\|_2^{1/2} \|\nabla u\|_2 + C \|\nabla u\|_2 \| \pp_z
T_2\|_2^{1/2} \|\nabla \pp_z T_2\|_2^{1/2} \|\tt\|_2^{1/2}
\|\nabla \tt\|_2^{1/2}.
\end{eqnarray*}
By Young's inequality, we get
\begin{eqnarray*}
&&\hskip-.68in \frac{d \|u\|_2^2}{dt} \leq C \left( \|\nabla
v_2\|_2^4+\|\nabla T_2\|_2^4+ \| \pp_z v_2\|_2^{2} \|\nabla \pp_z
v_2\|_2^{2} + \| \pp_z T_2\|_2^{2} \|\nabla \pp_z T_2\|_2^{2}
\right) \left( \|u\|_2^2 +\|\tt\|_2^2 \right).
\end{eqnarray*}
Thanks to Gronwall inequality, we obtain
\begin{eqnarray*}
&&\hskip-.68in \|u(t)\|_2^2 +\|\tt(t)\|_2^2
\leq \left( \|u(t=0)\|_2^2 +\|\tt(t=0)\|_2^2 \right) \times  \\
&&\hskip-.6in \exp \left\{ C  \int_0^t \left( \|\nabla v_2
(s)\|_2^4 +\|\nabla T_2(s)\|_2^4 + \| \pp_z v_2(s)\|_2^{2}
\|\nabla \pp_z v_2(s)\|_2^{2} +\| \pp_z T_2(s)\|_2^{2} \|\nabla
\pp_z T_2(s)\|_2^{2} \right) \; ds \right\}.
\end{eqnarray*}
Since $(v_2, T_2)$ is a strong solution,  we have
\begin{eqnarray*}
&&\hskip-.68in \|u(t)\|_2^2 +\|\tt(t)\|_2^2 \leq \left(
\|u(t=0)\|_2^2 +\|\tt(t=0)\|_2^2 \right) \, \exp\{ C  \left( K_V^2
t +K_t^2 t + K_z K_V  +K_t^{2} \right)  \}.
\end{eqnarray*}
The above inequality proves the continuous dependence of the
solutions on the initial data, and in particular, when
$u(t=0)=\tt(t=0)=0$, we have $u(t)=\tt(t)=0,$ for all $t\ge 0$.
Therefore, the strong solution is unique.

\end{proof}

\noindent
\section*{Acknowledgments}
We are thankful to the anonymous referee for the useful comments
and suggestions. This work was supported in part by the NSF grants No. DMS--0204794
and  DMS--0504619,  the MAOF Fellowship of the Israeli Council of
Higher Education, and by the USA Department of Energy, under
contract number W--7405--ENG--36 and ASCR Program in Applied
Mathematical Sciences.

\end{document}